\documentclass[10pt]{article}

\newcommand{\ddate}{March 21 2006}
\date{\ddate}

\usepackage{amssymb,amsmath}
\usepackage[matrix,arrow,curve]{xy}
\usepackage{epsfig}

\usepackage{pstricks}
\usepackage{epic}
\usepackage{eepic}
\newtheorem{thm}{Theorem}[section]

\newtheorem{pro}[thm]{Proposition}

\newtheorem{Theorem}[thm]{Theorem}

\newtheorem{Lemma}[thm]{Lemma}

\newtheorem{Proposition}[thm]{Proposition}

\newtheorem{Corollary}[thm]{Corollary}

\newtheorem{ccote}[thm]{\hskip -1.6mm}
\newcommand{\preu}{\noindent{\sc Proof: \ }}

\newcommand{\ppreu}[1]{\noindent {\sc Proof of #1: \ }}

\newcommand{\cqfd}{\unskip\kern 6pt\penalty 500
\raise -2pt\hbox{\vrule\vbox to10pt{\hrule width
 4pt\vfill\hrule}\vrule}\smallskip}

\newcommand{\proref}[1]{Proposition~\ref{#1}}

\newcommand{\lemref}[1]{Lemma~\ref{#1}}
\newcommand{\corref}[1]{Corollary~\ref{#1}}
\newcommand{\thref}[1]{Theorem~\ref{#1}}
\newcommand{\exref}[1]{Example~\ref{#1}}

\newcommand{\bbs}{{\mathbb{S}}}

\newcommand{\bbr}{{\mathbb{R}}}

\newcommand{\bbc}{{\mathbb{C}}}

\newcommand{\cala}{{\mathcal A}}

\newcommand{\calc}{{\mathcal C}}

\newcommand{\calh}{{\mathcal H}}

\newcommand{\calp}{{\mathcal P}}

\newcommand{\cals}{{\mathcal S}}

\newcommand{\pcirc}{\kern .7pt {\scriptstyle \circ} \kern 1pt}

\newcommand{\hfl}[1]{\buildrel{#1}\over{\longrightarrow}}
\newcounter{exo}

\newcommand{\mun}{{-1}}
\newcommand{\sk}[1]{\vskip #1 mm}

\newcommand{\llangle}[2]{\langle #1 ,#2 \rangle}

\newcommand{\onto}{\to\kern-7.5pt\to}
\newcommand{\donto}{\downarrow\kern -7.92pt\raisebox{-0.6ex}{$\downarrow$}}

\newcommand{\lso}{\mathfrak{so}}

\newcommand{\conf}{{\rm Conf\,}}
\newcommand{\moeb}[1]{{\rm M\ddot{o}b}({#1})}
\newcommand{\lmoeb}[1]{\mathfrak{m}\ddot{\mathfrak o}{\mathfrak b}({#1})}
\newcommand{\moe}{{\rm M\ddot{o}b}(d-1)}
\newcommand{\lmoe}{\mathfrak{m}\ddot{\mathfrak o}{\mathfrak b}(d-1)}
\newcommand{\Grad}{{\rm Grad\,}}
\newcommand{\grad}{{\rm grad\,}}

\newcommand{\ho}{{\rm horb}}
\newcommand{\spdim}{{\rm spdim}}
\newcommand{\sed}{{\rm sed}}

\newsavebox{\cer}

\title{Holonomy orbits of the snake charmer algorithm}
\author{Jean-Claude HAUSMANN and Eugenio RODRIGUEZ}
\begin{document}
\maketitle
%\tableofcontents

%\section{Introduction}\label{intro}

A {\it snake} (of length $L$) is a (continuous) piecewise
$\calc^1$-curve $\cals:[0,L]\to\bbr^d$,
parameterized by arc-length and whose ``tail'' is at the origin
($\cals(0)=0$). Charming a snake consists in having it move in such a way that
its ``snout'' $\cals(L)$ follows a chosen $\calc^1$-curve $\gamma(t)$.
The {\it snake charmer algorithm},  initiated in \cite{Ha2}
for polygonal snakes and developed in \cite{Ro} in the general case,
works as follows. The input is a pair $(\cals,\gamma)$, where:
\renewcommand{\labelenumi}{(\roman{enumi})} %%after "enumerate"
\begin{enumerate}
\item $\cals:[0,L]\to\bbr^d$ is a snake of length $L$,
\item $\gamma:[0,1]\to\bbr^d$ is $\calc^1$-curve with
$\gamma(0)=\cals(L)$.
\end{enumerate}
The output will then be a continuous $1$-parameter family
$\cals_t$ of snakes of length $L$ satisfying $\cals_0=\cals$
and $\cals_t(L)=\gamma(t)$. This algorithm, described in
Section~\ref{Salgo} below, is Ehresmannian in nature: the output
is a horizontal lifting for some
connection. A holonomy phenomenon for closed curves $\gamma$ then occurs:
having $\gamma(1)=\gamma(0)$ does not imply that $\cals_1=\cals_0$.
Given a snake $\cals$, one can input in the snake charmer algorithm
the pairs $(\cals,\gamma)$ for all possible $\calc^1$-loops $\gamma$ at $\cals(L)$.
The snakes $\cals_1$ obtained this way form the
{\it holonomy orbit} of $\cals$.
The purpose of this note
is to study these holonomy orbits, proving that, in good cases,
they are compact smooth manifolds diffeomorphic to real Stiefel manifolds.

The paper is organized as follows. After some preliminaries in
Section~\ref{Spreli}, we give in Section~\ref{Salgo} a survey
of the snake charmer algorithm and of its main properties
(more details are to be found in~\cite{Ro}). Section~\ref{Sorbits}
presents our new results about holonomy orbits with their proofs.
Section~\ref{Sexps} contains examples.

The first author is grateful to Giorgi Khimshiashvili
for a stimulating discussion during the conference.

%%%%%%%%%%%%%%%%%%%%%%%%%%%%%%%%%%%%%%%%%%%%%%%%%%%%%%%%%%%%%%%%%%%%%%%%%%%%%%
\section{Preliminaries}\label{Spreli}

\begin{ccote}\label{p-cont}\rm
Let $L$ be a positive real number and let
$\calp=\{0=s_0<s_1<\cdots<s_{N-1} < s_N=L\}$ be a finite partition of $[0,L]$.
If $(X,d)$ is a metric space, a map $z:[0,L]\to X$ is said
{\it piecewise continuous for $\calp$} if,
for every $i=0,\ldots,N-1$, the restriction of $z$ to the semi-open
interval $[s_i,s_{i+1})$ extends to a (unique) continuous map
$z_i$ defined on the closed interval $[s_i,s_{i+1}]$. In particular,
$z$ is continuous on the right and the discontinuities, only possible
at points of $\calp$, are just jumps.
We denote by $\calc^0_\calp([0,L],X)$ the space of maps from
$[0,L]$ to $X$ which are piecewise continuous for $\calp$; this space is
endowed with the uniform convergence distance
$$d(z_1,z_2)=\sup_{s\in [0,L]}\{z_1(s),z_2(s)\}\ .$$
When $\calp$ is empty, the map $z$ is just continuous and
$\calc^0_\emptyset([0,L],X)=\calc^0([0,L],X)$.
The map $z\mapsto (z_1,\dots,z_N)$ provides a homeomorphism
\begin{equation}\label{homeoprod}
\calc^0_\calp([0,L],X)\approx \prod_{i=0}^{N-1} \calc^0([s_i,s_{i+1}],X)\ .
\end{equation}
If $X$ is a Riemannian manifold, the space $\calc^0([s_i,s_{i+1}],X)$
naturally inherits a Banach manifold structure (see~\cite{Ee}),
so $\calc^0_\calp([0,L],X)$ is a Banach manifold using~\eqref{homeoprod}.
The tangent space $T_z\calc^0_\calp([0,L],X)$ to $\calc^0_\calp([0,L],X)$
at $z$ is the space of those $v\in\calc^0_\calp([0,L],TX)$
satisfying $p\pcirc v = z$, where $p:TX\to X$ is the natural projection.
\end{ccote}

\begin{ccote}\label{moeb}\rm
The unit sphere in $\bbr^d$ centered at the origin
is denoted by $\bbs^{d-1}$.
Let $\moe$ be the group of M\"obius transformations of
$\bbs^{d-1}$. It is a Lie group of dimension $d(d+1)/2$,
with $SO(d)$ as a compact maximal subgroup.
Its Lie algebra is denoted by $\lmoe$.
For $0\neq v\in\bbr^d$,
we define a $1$-parameter subgroup $\Gamma^v_t$ of $\moe$ by
$$
\Gamma^v_t=\varphi_v^{-1}\pcirc\rho^v_t\pcirc\varphi_v
$$
where
\begin{itemize}
\item $\varphi_v:\bbs^{d-1}\to\widehat\bbr^d$ is the stereographic projection
sending $v/|v|$ to $\infty$ and $-v/|v|$ to $0$;
\item $\rho^v_t:\widehat\bbr^d\to\widehat\bbr^d$ is the homothety $\rho^v_t(x)=e^{t|v|}x$.
\end{itemize}
Thus, $\Gamma^v_t$ is a purely hyperbolic flow with stable fixed point $v/|v|$
and unstable fixed point $-v/|v|$. We agree that $\Gamma^0_t = {\rm id}$.
Let $C_v\in\lmoe$ such that
$\Gamma^v_t=\exp(tC_v)$ (with $C_0=0$). The correspondence $v\mapsto C_v$
gives an injective linear map $\chi:\bbr^d\to\lmoe$; its image is a $d$-dimensional
vector subspace $\calh$ of $\lmoe$, supplementary to $\lso(d)$,
and $\calh$ generates $\lmoe$ as a Lie algebra. Let $\Delta^\calh$ be the right invariant distribution on $\moe$ which is equal to $\calh$ at the unit element.

The $1$-parameter subgroups $\Gamma^v_t$ may be used to built up a
diffeomorphism $\Psi:\bbr^d\times SO(d-1)\hfl{\approx}{}\moe$ defined by
$\Psi(v,\rho)=\Gamma^v_1\cdot\rho$.
\end{ccote}

%%%%%%%%%%%%%%%%%%%%%%%%%%%%%%%%%%%%%%%%%%%%%%%%%%%%%%%%%%%%%%%%%%%%%%%%%%%%%%
\section{The algorithm}\label{Salgo}
In this section, we give a survey of the snake charmer algorithm and
some of its properties. For details, see~\cite{Ro}.

\begin{ccote}\rm
Fix a positive real number $L$ and a finite set $\calp\subset [0,L]$
as in~\ref{p-cont}.
Let $\conf=\calc^0_\calp([0,L],\bbs^{d-1})$, with its Banach manifold structure
coming from the standard Riemannian structure on $\bbs^{d-1}$.
The inclusion of $\bbs^{d-1}$ into $\bbr^d$ makes $\conf$ a submanifold
of the Banach space $\prod_{i=0}^{N-1} \calc^0([s_i,s_{i+1}],\bbr^{d})$.

The space $\conf$ is the {\it space of configurations} for the snakes
of length $L$ which are continuous and ``piecewise $\calc^1$ for $\calp$''.
The snake $\cals_z$ associated to
$z\in\conf$ is the map $\cals_z(s)=\int_0^s\!z(\tau)d\tau$.
Taking its ``snout'' $S_z(L)$ provides a map $f:\conf\to\bbr^d$,
defined by
$$
f(z)=\int_0^L\! z(s)ds
$$
which is proven to be smooth. The image of $f$ is the closed ball
of radius $L$ centered at the origin.

The snakes corresponding to piecewise constant configurations,
$z(s) = z_i$ for $s \in [s_{i-1},s_i)$, are called {\it polygonal snakes}.
In this case, we see $f$ as a map from $(\bbs^{d-1})^N$ to $\bbr^d$ sending
$z = (z_1,\ldots,z_N)$ to $\sum_{i=1}^N (s_{i} - s_{i-1})z_i$. For more details
on this particular case see~\cite{Ha2}. If all the $(s_{i} - s_{i-1})$ are
equal, the snake is called an {\it isosceles polygonal snake}.

The critical points of $f$ are the {\it lined} configurations,
where $\{z(s)\mid s\in[0,L]\} \subset \{ \pm p \}$ for a point $p \in \bbs^{d-1}$
(for polygonal snakes, this is \cite[Theorem~3.1]{Ha1}).
The snake associated to such a configuration is then contained in the line through $p$.
The set of critical values is thus
a finite collection of $(d-1)$-spheres centered at the origin, which
depends on $\calp$.

Charming snakes will now be a path-lifting ability
for the map $f$: given an initial configuration $z \in \conf$
and a $\calc^1$-curve $\gamma : [0,1] \to \bbr^d$
such that $\gamma(0) = f(z)$, we are looking for a curve $t \mapsto z_t \in \conf$
such that $z_0 = z$ and $f(z_t) = \gamma(t)$.
In Ehresmann's spirit, we are looking for a connection
for the map $f$.
The tangent space $T_z\conf$ to $\conf$ at $z$ is the vector space
of those maps $v\in\calc^0_\calp([0,L],\bbr^d)$
such that $\llangle{v(s)}{z(s)}=0$ for all $s\in[0,L]$, where
$\llangle{}{}$ denotes the usual scalar product in $\bbr^d$.
We endow $T_z\conf$ with the scalar product
$\llangle{v}{w}=\int_0^L\!\llangle{v(s)}{w(s)}ds$. For each smooth map $\varphi:\bbr^d\to\bbr$,
one gets a vector field $\Grad(\varphi\pcirc f)$ on $\conf$ defined by
$$
\Grad_z(\varphi\pcirc f)(s) = \grad_{f(z)}\varphi  -
\llangle{z(s)}{\grad_{f(z)}\varphi}z(s) \ .
$$
This vector field plays the role of the gradient of $\varphi\pcirc f$, that is
$$
\llangle{\Grad_z(\varphi\pcirc f)}{v} = T_z(\varphi\pcirc f)(v)
$$
for each $v\in T_z\conf$
(as the metric induced by our scalar product
is not complete, gradients do not exist in general).
For $z\in\conf$, the set of all $\Grad_z(\varphi\pcirc f)$
for $\varphi\in\calc^1(\bbr^d,\bbr)$ is a vector subspace $\Delta_z$ of $T_z\conf$,
of dimension $d-1$ if $z$ is a lined configuration and $d$ otherwise.
The correspondence $z\mapsto\Delta_z$ is a distribution $\Delta$ (of non-constant dimension).
For a pair $(z,\gamma)\in\conf\times\calc^1([0,1],\bbr^d)$
such that $f(z)=\gamma(0)$, the snake charmer algorithm takes for $z_t$
the horizontal lifting of $\gamma$ for the connection $\Delta$.

As the map $f$ is not proper and the dimension of $\Delta_z$ is not constant,
the existence of horizontal liftings has to be established.
We use the $\calc^\infty$-action of the M\"obius group $\moe$ on $\conf$ by
post-composition: $g\cdot z = g\pcirc z$. For $z_0\in\conf$, let $\cala(z_0)$
be the subspace of those $z\in \conf$ which can be joined to $z_0$ by
a succession of $\Delta$-horizontal curves.
One of the main results
(\cite[Theorem~2.19]{Ro}, proven in~\cite{Ha2}
for isosceles polygonal snakes)
says that $\cala(z_0)$ coincides with the orbit of $z_0$
under the action of $\moe$:
\begin{Theorem}\label{thpri}
For all $z_0\in\conf$, one has $\cala(z_0)=\moe\cdot z_0$.
\end{Theorem}
The proof of \thref{thpri} uses the following fact about the action of
the $1$-parameter subgroup $\Gamma^v_t$ of $\moe$ introduced in~\ref{moeb}:
\begin{Lemma}\label{Lflot}
The flow $\phi_t(z)=\Gamma^v_t\cdot z$ on $\conf$ is the gradient flow
of the map $z\mapsto\llangle{f(z)}{v}$.
\end{Lemma}
The flow $\phi_t$ on $\conf$ is therefore $\Delta$-horizontal.
If $z\in\conf$, denote by $\beta_z:\moe\to\conf$ the smooth map
$\beta_z(g)=g\cdot z$. \lemref{Lflot} implies that the map $\beta_z$
sends the distribution $\Delta^\calh$ onto the distribution $\Delta$.
If $z$ is not a lined configuration, $T_g\beta_z :\Delta^\calh_g\to\Delta_{g\cdot z}$
is an isomorphism; if $z$ is lined, there is a $1$-dimensional kernel
(for instance, if $z(s)\in\{\pm p\}$, then $\Gamma^p_t\cdot z = z$ for all $t$).
\end{ccote}

\begin{ccote}\label{ode-par} The differential equation.\ \rm
As a consequence of \thref{thpri}, a $\Delta$-horizontal curve
$z_t$ starting at $z_0$ may be written as $g(t) \cdot z_0$ where $g(t) \in \moe$
and $g(0) = {\rm id}$. For a given $\calc^1$-curve $\gamma(t) \in \bbr^d$,
the $\Delta$-horizontal lifting $g(t) \cdot z_0$ of $\gamma(t)$ is
obtained by taking for $g(t)$ the solution of an ordinary differential equation
in the M\"obius group $\moe$ that we describe here below
(more details can be found in~\cite[Chapter 3]{Ro}).

For $z\in\conf$, we define the $(d \times d)$-matrix $M(z)$ by
$$
M(z) = \begin{pmatrix} L & & 0 \\
& \ddots & \\
0 &  & L
\end{pmatrix} -
\begin{pmatrix}
\int_0^L z_1(s)z_1(s)ds & \cdots & \int_0^L z_1(s)z_d(s)ds \\
\vdots & & \vdots \\
\int_0^L z_1(s)z_d(s)ds & \cdots & \int_0^L z_d(s)z_d(s)ds
\end{pmatrix} \ .
$$
(Observe that the second term is the Gram matrix of
the vectors $(z_1,\ldots,z_d)$ for the scalar product).
It can be proved that $M(z)$ is invertible if and only if $z$ is not lined.
Let $F = F_{z_0}:\moe\to\bbr^d$ be the composition $F=f\pcirc\beta_{z_0}$.
Let $g\in\moe$. Using the linear injective map
$\chi$ of~\ref{moeb}, with image $\calh$, and the right translation
diffeomorphism $R_g$ in $\moe$, we get a linear map
\begin{equation}\label{linmapM}
\begin{array}{c}{\xymatrix@C-3pt@M+2pt@R-4pt{%
\bbr^d\ar[r]^(0.50){\chi}_(0.50){\approx} &
\calh \ar[r]^(0.40){T_eR_g}_(0.40){\approx} &
\Delta^\calh_g \ar[r]^(0.30){T_{g}F} &
T_{F(g)}\bbr^d \approx \bbr^d
}}\end{array}
\end{equation}
It turns out that the matrix of the linear map~\eqref{linmapM} is
$M(g\cdot z_0)$. If $g\cdot z_0$ is not a lined configuration,
$T_gF$ is bijective and $M(g\cdot z_0)$ is invertible.
One way to insure that $g\cdot z_0$ is not lined for all
$g\in\moe$ is to assume that $z_0$ takes at least $3$ distinct values.
The following result is proven
in \cite[Prop.\,3.10]{Ro}:
\begin{Proposition}\label{ode_mobeq_Pro}
Let $z_0\in\conf$ be a configuration that takes at least $3$ distinct values.
Let $\gamma:[0,1]\to\bbr^d$ be a $\calc^1$-curve with $\gamma(0)=f(z)$.
Suppose that
the $\calc^1$-curve $g:[0,1]\to\moe$ satisfies the following differential
equation:
\begin{equation}\label{ode_mob}
\dot{g}(t) = R_{g(t)}\pcirc\chi\Big(
M^{-1}(g(t) \cdot z_0)\,\dot{\gamma}(t)\Big) \ , \quad g(0) = {\rm id} \ .
\end{equation}
Then, the curve $g(t)\cdot z_0$ is the unique $\Delta$-horizontal lifting of
$\gamma$ starting at $z_0$.
\end{Proposition}

\iffalse
The idea of the proof of \proref{ode_mobeq_Pro} is the following.
Denote by $X^t$ the map
$g \mapsto R_{g(t)}\pcirc\chi\big(M^{-1}(g(t) \cdot z_0)\,\dot{\gamma}(t)\big)$.
This is
that a time dependent vector field on $\moe$ and Equation~\eqref{ode_mob}
becomes $\dot{g}(t) = X^t(g(t))$.
This vector field is $\calc^\infty$ in the variable $g$ (since the
action of $\moe$ on $\conf$ is $\calc^\infty$) but only continuous
in $t$. By direct calculation it can be shown that the derivative
of $X^t$ (in the variable $g$) is continuous in $t$ and thus classical results
concerning existence and uniqueness of solutions
of ordinary differential equations may be applied.
In~\ref{contSol}, we give more details about how the solution $g(t)$
depends on $z_0$ and $\gamma$.
\fi

As the M\"obius group is not compact and the map $F$
is not proper, Differential equation~\eqref{ode_mob}
may not have a solution for all $t\in [0,1]$. The notion of ``sedentariness'',
described in~\ref{sedentariness}, is the main tool to study the global
existence of $\Delta$-horizontal liftings.
\end{ccote}

\begin{ccote}\label{odedeq2} The cases $d=2$ or $d=3$.\ \rm
For planar snakes ($d=2$), one can use that $\moeb{1}$ is isomorphic to
$PSU(1,1) = SU(1,1)/\{\pm {\rm id} \}$,
where $SU(1,1) = \{ \left(\begin{smallmatrix} a & b \\ \bar{b} & \bar{a} \end{smallmatrix}\right)
\mid a,b \in \bbc \text{ and } |a|^2 - |b|^2 = 1 \}$.
The isomorphism between $\moeb{1}$ and $PSU(1,1)$ comes from the
$SU(1,1)$-action on $\bbc$ by homographic transformations,
which preserves the unit circle.
The group $SU(1,1)$ being thus a $2$-fold covering of $\moeb{1}$, the curve
$g(t)\in\moeb{1}$ of \proref{ode_mobeq_Pro} admits a unique lifting in
$\tilde g(t)\in SU(1,1)$ with $\tilde g(0)={\rm id}$. The output
of the snake charmer algorithm will then be of the form $z_t=\tilde g(t)\cdot z_0$,
using the action of $SU(1,1)$ on $\conf(1)$. Working in the matrix group
$SU(1,1)$ is of course favorable for numerical computations.

The Lie algebra $\mathfrak{psu}(1,1) = \mathfrak{su}(1,1)$
consists in matrices of the form $\left(\begin{smallmatrix} iu & b \\ \bar{b} & -iu \end{smallmatrix}\right)$ with $u\in\bbr$ and $b\in\bbc$. 
Under the isomorphism $\mathfrak{su}(1,1)\approx\lmoeb{1}$, the element
$\chi(v)\in\calh\subset\lmoeb{1}$ corresponds to the matrix
$\left(\begin{smallmatrix} 0 & v \\ \bar{v} & 0 \end{smallmatrix}\right) 
\in\mathfrak{su}(1,1)$. 
These matrices form a $2$-dimensional vector space 
$\tilde\calh\subset\mathfrak{su}(1,1)$
giving rise to a right invariant distribution $\Delta^{\tilde\calh}$
on $SU(1,1)$. By \proref{ode_mobeq_Pro}, the curve $\tilde g(t)\cdot z_0\in\conf(1)$
is horizontal if $\tilde g(t)$ is $\Delta^{\tilde\calh}$-horizontal.
In this language, Equation~\eqref{ode_mob} becomes
\begin{equation}\label{ode_mob_d2}
\dot{\tilde g}(t) = \begin{pmatrix}
0 & v/2 \\ \bar{v}/2 & 0
\end{pmatrix} \tilde g(t), \quad\quad g(0) = \begin{pmatrix}
1 & 0 \\ 0 & 1
\end{pmatrix}
\end{equation}
where $v = v_1(g,t) + i v_2(g,t)$
is obtained by
$$
\begin{pmatrix}
v_1(g(t),t) \\
v_2(g(t),t)
\end{pmatrix} = M(\tilde g(t) \cdot z_0)
\begin{pmatrix}
\dot{\gamma}_1(t) \\
\dot{\gamma}_2(t)
\end{pmatrix}
$$
For more details, see~\cite[Chapter 5]{Ro}. An analogous 
(and formally similar) matrix 
approach is available for spatial snakes ($d=3$), using
the isomorphism $PSL(2,\bbc)\approx\moeb{2}$;
see~\cite[Prop.~5.6]{Ro}.
\end{ccote}

\begin{ccote}\label{sedentariness}Sedentariness.\ \rm
Denote by $\mu$ the Lebesgue measure on $\bbr$. Let $z\in\conf$. The
{\it sedentariness} $\sed(z)\in[0,L]$ of $z$ is defined by
$\sed(z)=\max_{p\in \bbs^{d-1}} \mu(z^\mun(p))$. This maximum exists since
the set $\{p\in \bbs^{d-1}\mid \mu(z^\mun(p))>r\}$ is finite for all
$r>0$. By \thref{thpri}, $\sed(z)$ is an invariant of $\cala(z)$.
Observe that, if $\sed(z)\neq L/2$, then $z$ takes at least 3 distinct values.
The sedentariness of $z$ is used in~\cite[Section~3.3]{Ro}
to get global existence results
for horizontal liftings of a path starting at $f(z)$. The one we need
is the following

\begin{Proposition}\label{threlev}
Let $z\in\conf$ and let $\gamma:[0,1]\to\bbr^d$ be a $\calc^1$-path
with $\gamma(0)=f(z)$. Suppose that $\gamma([0,1])$ is
contained in the open ball centered at $0$ of radius $L-2\,\sed(z)$.
Then there exists a (unique) $\Delta$-horizontal lifting
$\tilde\gamma:[0,1]\to\conf$ for $\gamma$ with $\tilde\gamma(0)=z$.
\end{Proposition}

As the sedentariness of a configuration is preserved along horizontal
curves, \proref{threlev} is also true for continuous piecewise
$\calc^1$-paths.
A configuration $z$ is called {\it nomadic} if $\sed(z)=0$. \proref{threlev}
guarantees that, if $z$ is nomadic, any $\calc^1$-path
starting at $f(z)$ admits a horizontal lifting, provided its image stays
in the open ball of radius $L$. More general results may be found
in~\cite[Section~3.3]{Ro}.
\end{ccote}

\begin{ccote}\label{contSol} Continuity of the algorithm.\ \rm
We now describe how the snake charmer algorithm behaves as we vary the initial
configuration $z_0$ and the $\calc^1$-curve $\gamma$. Our goal is not
to present the most general statements but only those needed in
Section~\ref{Sorbits}.

Let $z_0 \in \conf$, $b = f(z_0)$ and let $\conf^b=f^\mun(b)$.
Set $\sigma = \sed(z_0)$ and suppose that $|f(z_0)| < L - 2\sigma$.
Denote by $B_{L-2\sigma}(0)$ the open ball in $\bbr^d$ of radius $L-2\sigma$
centered at the origin. Consider a $\calc^1$-curve
$\gamma : [0,1] \to B_{L-2\sigma}(0)$ with $\gamma(0) = b$.
Define $E_{\sigma,b} = \{ z \in \conf^b \mid \sed(z) \leq \sigma \}$;
it is a metric space with the induced metric from $\conf$.
Consider the map
$$
X : \moe \times E_{\sigma,b} \times [0,1] \longmapsto T\moe,
$$
defined by
$$
X(g,z,t)=R_{g(t)}\pcirc\chi\Big(M^{-1}(g(t) \cdot z_0)\,\dot{\gamma}(t)\Big)\ .
$$
This is a vector field on $\moe$ depending on the time $t$ and
on the parameter $z\in E_{\sigma,b}$.
Differential equation~\eqref{ode_mob} becomes
\begin{equation}\label{ode_mob_par}
\dot{g}(t) = X(g(t),z,t) \ , \quad g(0) = {\rm id} \ .
\end{equation}
For any given $z \in E_{\sigma,b}$, the solution $g_z$ of Equation~\eqref{ode_mob_par}
exists for all $t \in [0,1]$ by \proref{threlev}.
Since $X$ is continuous and its derivative in $g$ is a continuous map on the variables $z$ and $t$,
classical results on the dependence of parameters for the solution
(see for example~\cite[4.3.11]{Sc}) imply that the map $(z,t) \mapsto g_z(t)$
is continuous. This yields:

\begin{pro}[Continuity in $z$]\label{continuite}
Let $\gamma : [0,1] \to B_{L-2\sigma}(0)$ be a $\calc^1$-curve
such that $\gamma(0) = b$. For any $z \in E_{\sigma,b}$, denote by
$\tilde{\gamma}_z(t) = g_z(t) \cdot z$ the
$\Delta$-horizontal lift of $\gamma$ starting at $z$.
The map $E_{\sigma,b} \times [0,1] \to \conf$
that sends $(z,t)$ to $\tilde{\gamma}_z(t)$ is continuous.
\end{pro}

Let $M \subset E_{\sigma,b}$ be a smooth submanifold of $\conf$ and
suppose that $\gamma$ is of class $\calc^2$. The map $X$ restricted to $M$ is therefore
$\calc^1$ (in all variables) and using once again~\cite[4.3.11]{Sc} we get:

\begin{pro}[Differentiability in $z$]\label{differentiability}
Let $\gamma : [0,1] \to B_{L-2\sigma}(0)$ be a $\calc^2$-curve
such that $\gamma(0) = b$. For any $z \in M$, denote by
$\tilde{\gamma}_z(t) = g_z(t) \cdot z$ the
$\Delta$-horizontal lift of $\gamma$ starting at $z$.
The map $M \times [0,1] \to \conf$
that sends $(z,t)$ to $\tilde{\gamma}_z(t)$ is of class $\calc^1$.
\end{pro}

We have an analog result when we fix the configuration
$z_0$ and vary $\gamma$. Consider the set $C_{\sigma,b}
= \{ \gamma \in \calc^1([0,1],\bbr^d) \mid \gamma(0) = b, \; \gamma([0,1]) \subset
B_{L-2\sigma}(0) \}$ that is an open subset of the affine space of $\calc^1$-paths
starting at $b$.

\begin{pro}[Continuity in $\gamma$]\label{continuite_gamma}
Let $z_0 \in \conf$, $b = f(z_0)$ and $\sigma = \sed(z_0)$.
Suppose that $|f(z_0)| < L - 2\sigma$. For any $\gamma \in C_{\sigma,b}$,
denote by $\tilde{\gamma}_z(t)$ the
$\Delta$-horizontal lift of $\gamma$ starting at $z_0$.
The map $C_{\sigma,b} \times [0,1] \to \conf$
that sends $(\gamma,t)$ to $\tilde{\gamma}_z(t)$ is continuous.
\end{pro}
\end{ccote}

\begin{ccote}\label{bival} Bivalued configurations.\ \rm
Let $z_0$ be a bivalued configuration, that is $z_0([0,L])=\{p_0,q_0\}$,
with $p_0\neq q_0$. Let $L_p$ and $L_q$ be the Lebesgue measures of
$z_0^\mun(p_0)$ and $z_0^\mun(q_0)$ (these preimages are finite unions of
intervals). By \thref{thpri},
a horizontal curve $z_t$ starting at $z_0$ will stay bivalued:
$z_t([0,L])=\{p(t),q(t)\}$, with $p(0)=p_0$, $q(0)=q_0$ and $p(t)\neq q(t)$.
Also, one has $z_t^\mun(p(t))=z_0^\mun(p_0)$, $z_t^\mun(q(t))=z_0^\mun(q_0)$
and therefore $z_t$ is determined by the pair $(p(t),q(t))$.
Hence, $\cala(z_0)$ is contained in a compact submanifold $W$ of $\conf$
naturally parameterized by $\bbs^{d-1}\times \bbs^{d-1}$.
The restriction
of $f$ to $W$ gives a $\calc^\infty$-map $\hat f:W\to\bbr^d$ which takes
the explicit form $\hat f(p,q)=L_pp+L_qq$.
The definition of the distribution $\Delta$ implies that
$\Delta_{(p,q)}\subset T_{(p,q)}W$ when $(p,q)\in W$.

Consider the open sets $W^0\subset W$ and
$U^0\subset \bbr^d$ defined by
$$
W^0:=\{(p,q)\in W\mid p\neq \pm q\} \;,\quad
U^0:=\{x\in\bbr^d\mid |x|< L \hbox{ and } |x|\neq L_p-L_q \}
$$
As $W^0$ contains no lined configurations,
the map $\hat f$ restricts to a submersion $\hat f^0=f^0:W^0\to U^0$ for which
$\Delta$ is an ordinary connection. As $f^0$
extends to $\hat f:W\to\bbr^d$ and as $W$ is compact, the map
$f^0$ is proper. The original Ehresmann's construction of
horizontal liftings \cite{Eh} then apply: if $z_0\in W^0$ and
$\gamma:[0,1]\to U^0$ is a $\calc^1$-curve
$\gamma(0)=z_0(L)$, then $\gamma$ admits a
$\Delta$-horizontal lifting $\tilde\gamma:[0,1]\to W^0$ of $f^0$
with $\tilde\gamma(0)=z_0$.

Let us specialize to planar snakes ($d=2$). Each fiber of $f^0$
then consists of two points; as $f^0$ is proper, it is thus
a $2$-fold cover of $U^0$. Therefore, $\tilde\gamma(t)$ is
determined by $f(\tilde\gamma(t))=\gamma(t)$. We deduce that
{\it the unique lifting of $\gamma$ into $\conf$ is horizontal}.

If, for $t=t_0$, $\gamma(t)$ crosses the sphere of radius $L_p-L_q$,
it may happen that $\tilde\gamma(t)$ tends to a lined configuration
when $t\to t_0$ (see \exref{exnoconnhorb} below). To understand
when the unique lifting $\tilde\gamma(t)$ is horizontal at $t_0$,
we must study horizontal liftings around a lined configuration,
which, after changing notations, we call again $z_0$, corresponding to
$(p_0,q_0)\in W$ with $p_0=-q_0$.
The vector space $\Delta_{z_0}$ is then of dimension $1$ and it turns
out that $T_{z_0} f(\Delta_{z_0})$ is the line orthogonal to $p_0$,
see \cite[Remark 1.17]{Ro}.
Let $\gamma:[0,1]\to\bbr^d$ be a $\calc^1$-curve with
$\gamma(0)=z_0(L)$ and $\dot\gamma(0)\neq 0$.
A necessary condition for $\gamma$ to admit a
horizontal lifting is then $\langle\dot\gamma(0),p_0\rangle=0$.
If we orient the plane with the basis $(p_0,\dot\gamma(0))$,
the curve $\gamma$ has a signed curvature $\kappa(0)$ at $t=0$.
In~\cite[Prop. 3.18]{Ro}, it is proven that, around $t=0$, the unique lifting
of $\gamma$ into $\conf$ is horizontal if and only if
$\dot\gamma(0)$ is orthogonal to $p_0$ and
\begin{equation}\label{Hairereq}
\kappa(0)=\frac{L_{p_0}-L_{q_0}}{L^2} .
\end{equation}
Condition~\eqref{Hairereq} has been detected by the numerician Ernst Hairer.
Such a second-order condition for the existence of a horizontal lifting
for a non-constant rank distribution is worth being studied.
As far as we know, no such a phenomenon is
mentioned in the literature.
\end{ccote}

\begin{ccote}\label{fin}Miscellaneous.\ \rm
We finish this section by listing a few more properties of the
snake charmer algorithm.
Let $(\cals,\gamma)$ be an input for the algorithm,
with $\cals$ a snake of length $L$.
Let $\cals_t:[0,L]\to\bbr^d$ be the output.
The following two results follow directly from \thref{thpri}.

\begin{Proposition}[Regularity]\label{Pregu}
Suppose that, on some open subset $U\subset [0,L]$,
the snake $\cals$ is of class $\calc^k$, $k\in\{1,2,\dots,\infty\}$.
Then the same holds true for the snake $\cals_t$ for all $t\in [0,1]$.
\end{Proposition}

\begin{Proposition}[Periodicity]\label{Pperio}
Suppose that there exists $T\in\bbr$ such that $z(s)=z(s+T)$ for
all
$s$ such that $s$ and $s+T$ belong to $[0,L]$. Then,
$z_t(s)=z_t(s+T)$ for all $t\in [0,1]$.
\end{Proposition}

The following proposition follows either from \thref{thpri} or simply
from the fact that $\cals_t$ is a horizontal lifting.

\begin{Proposition}[Reparameterization]\label{Prepara}
Let $\varphi:[0,1]\to [0,1]$ be an orientation preserving
$\calc^1$-diffeomorphism.
Then, the deformation of $\cals$ following the curve $\gamma(\varphi(u))$
is $\cals_{\varphi(u)}$
\end{Proposition}

Finally, by construction, the distribution $\Delta$ is orthogonal
to fibers of $f:\conf\to\bbr^d$.
This can be rephrased in the following

\begin{Proposition}\label{encin}
$z_t$ is the unique lifting of $\gamma$
which, for all $t$,  minimizes the infinitesimal kinetic energy
of the hodograph, that is
$\frac{1}{2}\int_0^L\!|\frac{d}{dt}z_t(s)|^2\, ds$.
\end{Proposition}

Again, the existence of such minimizer $\cals_t$ is only guaranteed
by an analysis like in~\ref{sedentariness}.
\end{ccote}

%%%%%%%%%%%%%%%%%%%%%%%%%%%%%%%%%%%%%%%%%%%%%%%%%%%%%%%%%%%%%%%%%%%%%%%%%%%%%%
\section{Holonomy orbits}\label{Sorbits}

Let $b\in\bbr^d$. Define $\conf^b=f^\mun(b)$, the space of all configurations
with associated snake ending at $b$. Define the {\it holonomy orbit}
$\ho(z_0)\subset \conf^b$ of $z_0\in\conf^b$ by $\ho(z_0)=\cala(z_0)\cap\conf^b$;
it is thus the subspace of those $z\in\conf$
which are the result of the holonomy of the snake charmer algorithm for a pair
$(z_0,\gamma)$ with $\gamma$ a piecewise $\calc^1$-loop at $b$.
Even if $z_0$ is not lined, the point $b$ might not be a regular value of $f$,
so $f^\mun(b)$ might not be a submanifold of $\conf$. But \lemref{Thorb-lem}
below tells us that this is the case for $\ho(z_0)$. For $z\in\conf$,
define the {\it spherical dimension} $\spdim(z)$ of $z$ to be the minimal
dimension of a sub-sphere of $\bbs^{d-1}$ containing the set $z([0,L])$.
By \thref{thpri}, $\spdim(z)$ is an invariant of $\cala(z)$.
Notice that $\spdim(z)=0$ if and only if the configuration $z$
takes one or two values.

\begin{Lemma}\label{Thorb-lem}
Let $z_0\in\conf$. Let $k=\spdim(z_0)$ and  suppose that $k>0$.
Then $\ho(z_0)$ is a smooth submanifold of $\conf$ of
dimension $\sum_{i=1}^{k+1}(d-i)$.
\end{Lemma}
Generically, $\spdim(z_0)=d-1$, so $\ho(z_0)$ is of dimension $d(d-1)/2$.
The case $\spdim(z_0)=0$ is not covered by \lemref{Thorb-lem}.
It contains the monovalued case ($z$ being constant)
where $\ho(z_0)=f^\mun(f(z_0))=\{z_0\}$. The other case,
formed by the bivalued configurations, is interesting and is treated
in~\exref{exnoconnhorb} and \proref{ho_bival}.

\sk{3}
\ppreu{\lemref{Thorb-lem}}
Let $\beta:\moe\to\conf$ be the map $\beta(g)=g\cdot z_0$ and let
$F:\moe\to\bbr^d$ be the composition $F=f\pcirc\beta$.
Let $K=F^\mun(f(z_0))$. By \thref{thpri}, one has $\ho(z_0)=\beta(K)$.
The map $F$ is smooth.

Since $\spdim(z_0)>0$, the configuration $z_0$ takes at least three
values in~$\bbs^{d-1}$ and so does $g\cdot z_0$ for all $g\in\moe$. Therefore,
$\cala(z_0)$ contains no lined configurations and the tangent map
$T_gF : T_g\moe\to T_{F(g)}\bbr^d \approx \bbr^d$ is surjective for all
$g\in\moe$ (restricted to $\Delta^\calh_g$, it is an isomorphism).
Hence, $F$ is a submersion. Thus, $K$ is a smooth submanifold
of $\moe$ of dimension $\dim\moe - d = d(d+1)/2-d = d(d-1)/2$.

The manifold $K$ contains the stabilizer $A$ of $z_0$ and
$K\cdot A =K$. Let $V^{k}$ be the smallest sub-sphere of $\bbs^{d-1}$
containing the image of $z_0$. The group $A$ is then the stabilizer
of the points of $V^{k}$; since $k>0$, $A$ is conjugate in $\moe$ to
$SO(d-k-1)$. Therefore, the quotient space $K/A$ is of dimension
\begin{equation}\label{Thorb-lem-eq}
\frac{d(d-1)}{2} - \frac{(d-k-1)(d-k-2)}{2} =
\sum_{i=1}^{d-1}i - \sum_{j=1}^{d-k-2}\!\!j = \sum_{i=1}^{k+1}(d-i) \ .
\end{equation}
By \cite[Proposition 2.33]{Ro}, the map $\beta$ induces
an embedding of $\moe/A$ into $\conf$, hence an
embedding of $K/A$ into $\conf$ with image $\ho(z_0)$.
This proves \lemref{Thorb-lem}. \cqfd

In general, $\ho(z)$ is not closed; see examples
in \cite[pp.\,112--113]{Ha2}\footnote{In~\cite{Ha2},
one studies the set $\calh_b(U) \cdot z_0$: if $b = f(z_0)$
and $U$ is an open neighborhood of $b$,
$\calh_b(U)$ denotes the subgroup of diffeomorphisms of $\conf^b$
obtained as the holonomy of a piecewise $\calc^1$-loop at $b$ that
stays in $U$. Notice that $\calh_b(U) \cdot z_0 \subset \ho(z_0)$
and that $\calh_b(U) \cdot z_0$ contains the connected component
of $z_0$ of $\ho(z_0)$.}.
But this is the case if $f(z)$ is near the origin:

\begin{Theorem}\label{Thorb}
Let $z_0\in\conf$ such that $|f(z_0)|<L-2\,\sed(z_0)$. Let $k=\spdim(z_0)$
and suppose that $k>0$.
Then, each connected component of $\ho(z_0)$ is a
closed smooth submanifold of $\conf$ which is
diffeomorphic to the homogeneous space $SO(d)/SO(d-k-1)$.
\end{Theorem}

Observe that $SO(d)/SO(d-k-1)$ is the Stiefel manifold of orthonormal
$(k+1)$-frames in $\bbr^d$.

\sk{2}
\preu Suppose first that $f(z_0)=0$.
Let $z\in\ho(z_0)$.
For each $g\in SO(d)$, one has $f(g\cdot z)=g\cdot f(z)=0$, hence
$SO(d)\cdot z\subset \ho(z_0) = \ho(z)$.
As in the proof of \lemref{Thorb-lem}, let $\beta_z:\moe\to\conf$
be the map $\beta_z(g)=g\cdot z$. As $\spdim(z)=\spdim(z_0)=k>0$,
the stabilizer of $z$ is conjugate in $SO(d)$ to $SO(d-k-1)$.
Therefore, $\beta_z(SO(d))$ is a connected, compact submanifold of $\conf$
diffeomorphic to $SO(d)/SO(d-k-1)$.
By \lemref{Thorb-lem}, $\beta_z(SO(d))$ must be the connected component
of $z$ in $\ho(z_0)$. This proves \thref{Thorb} when $f(z_0)=0$.

In the general case, let $b=f(z_0)$ and let
$\gamma:[0,1]\to \bbr^d$ be the
linear parameterization $\gamma(t)=tb$ of the segment joining $0$ to $b$.
If $\delta:[0,1]\to X$ is a map, one denotes by
$\delta^-:[0,1]\to X$ the map $\delta^-(t)=\delta(1-t)$.
Let $z\in\conf^b$ such that $\sed(z)\leq\sed(z_0)$.
By \proref{threlev}, the path $\gamma$ admits a $\Delta$-horizontal
lifting $\tilde\gamma_z$ in $\conf$, starting at $z$.
The $\Delta$-parallel transport
$\tau(z)=\tilde{\gamma}_z(1)$ is then defined. Of course,
$\gamma^-$ also admits a horizontal lifting
$\tilde\gamma^-_{\tau(z)}$ which is equal to $(\gamma_z)^-$.
By \proref{continuite}, the parallel transport gives rise to a homeomorphism
$$
\tau : \{z\in\conf^b\mid \sed(z)\leq\sed(z_0)\}
\stackrel{\approx}{\longrightarrow} \{z\in\conf^0\mid \sed(z)\leq\sed(z_0)\} \ .
$$
One has $\tau(\ho(z_0))=\ho(\tau(z_0))$. As $\ho(z_0)$ is a smooth submanifold
of $\conf$, $\tau_{|\ho(z_0)}$ and $\tau^\mun_{|\tau(\ho(z_0))}$
are smooth by \proref{differentiability}.
Therefore $\tau$ induces a diffeomorphism from $\ho(z_0)$
onto $\ho(\tau(z_0))$. This terminates the proof of \thref{Thorb}.
\cqfd

In the particular case of planar snakes ($d=2$), we get:

\begin{Corollary}\label{Thorb-planar}
Let $z_0\in\conf(1)$ such that $|f(z_0)|<L-2\,\sed(z_0)$ and
\\ $\spdim(z_0)=1$. Then $\ho(z_0)$ is a disjoint union of circles.
\end{Corollary}

We do not know, under the hypotheses of \thref{Thorb}
or \corref{Thorb-planar}, whether the manifold
$\ho(z_0)$ is connected. This is not true if $|f(z_0)|>L-2\,\sed(z_0)$,
see Example~\ref{exnoconnhorb}.
It is however the case if $\sed(z_0)=0$:

\begin{Proposition}\label{Thorb-nomadic}
Let $z_0\in\conf$ be a nomadic configuration. Let $k=\spdim(z_0)$.
Then, $\ho(z_0)$ is a closed smooth submanifold of $\conf$ which is
diffeomorphic to the homogeneous space $SO(d)/SO(d-k-1)$.
\end{Proposition}

To prove \proref{Thorb-nomadic}, we need the following lemma
whose proof is postponed till the end of this section.

\begin{Lemma}\label{cinfloops}
Let $z_0\in\conf$ and let $z\in\ho(z_0)$. Then, there
exists a loop $\gamma:[0,1]\to\bbr^d$ at $f(z_0)$ which is of class $\calc^\infty$
and which admits a horizontal lifting joining $z_0$ to $z$.
\end{Lemma}

\ppreu{\proref{Thorb-nomadic}}
A nomadic configuration is not lined. Therefore, $k>0$
and the condition $|f(z_0)| < L-2\,\sed(z_0)=L$
is automatic. By \thref{Thorb}, it is then enough
to prove that $\ho(z_0)$ is connected.

Let $z\in\ho(z_0)$. By \lemref{cinfloops}, there
exists a loop $\gamma:[0,1]\to\bbr^d$ at $f(z_0)$, of class $\calc^\infty$,
admitting a horizontal lifting joining $z_0$ to $z$.
For $s\in [0,1]$, define
$\gamma_s(t) = (1-s)\gamma(0) + s\gamma(t)$.
%$\gamma_s(t) = \gamma(0) + (\gamma(t) - \gamma(0))s$.
The map $s\mapsto\gamma_s$ is a
homotopy of $\calc^\infty$-loops at $f(z_0)$, from $\gamma$ to the
constant loop. By \proref{threlev}, each loop $\gamma_s$ admits a horizontal lifting
$\tilde\gamma_s$ starting at $z_0$. By \proref{continuite_gamma}, the curve
$s\mapsto\tilde\gamma_s(1)$ is continuous, producing a path in $\ho(z_0)$
from $z$ to $z_0$. This proves that $\ho(z_0)$ is connected, provided we prove
\lemref{cinfloops}.
\cqfd

\ppreu{\lemref{cinfloops}}
By \thref{thpri}, there exists $g\in\moe$ with $z=g\cdot z_0$.
The M\"obius group $\moe$ is arc-wise connected using only curves that are piecewise
trajectories of flows of the type $\Gamma_t^v$ (see~\cite[p.\,102]{Ha2}).
Therefore, there exists $(v_i,\lambda_i)\in\bbr^d\times\bbr$,
$i=1,\dots,r$, such that $g=\Gamma^{v_r}_{\lambda_r}\cdots\Gamma^{v_1}_{\lambda_1}$.
For $i=1,\dots,r$, let $\phi_i:[0,1]\to [0,1]$ be a $\calc^\infty$-function
such that $\phi_i([0,\frac{i-1}{r}])=\{0\}$ and $\phi_i([\frac{i}{r},1])=\{1\}$.
The map $g(t):[0,1]\to\moe$ given by
$$
g(t)=\Gamma^{v_r}_{\phi_r(t)\lambda_r}\cdots\Gamma^{v_1}_{\phi_1(t)\lambda_1}
$$
is a $\calc^\infty$-curve in $\moe$, joining the unit element to $g$.
By \lemref{Lflot}, the curve $g(t)$ is $\Delta^\calh$-horizontal.
Then $g(t) \cdot z_0$ is the horizontal lifting of the $\calc^\infty$-path
$\gamma(t)=f(g(t) \cdot z_0)$, which is a loop of class $\calc^\infty$
at $f(z_0)$. \cqfd

%%%%%%%%%%%%%%%%%%%%%%%%%%%%%%%%%%%%%%%%%%%%%%%%%%%%%%%%%%%%%%%%%%%%%%%%%%%%%%
\section{Examples}\label{Sexps}

All our examples consist of planar snakes ($d=2$). 
We identify $\bbr^2$ with the complex plane $\bbc$. Configurations
$z \in \conf(1)$ are expressed under the form of $s \mapsto e^{i\theta(s)}$
with $\theta(s) \in \bbr$.

\begin{ccote}\rm
In our first example, the snake is a half circle with
configuration $z:[0,\pi]\to \bbs^1$ given by $z(s)=ie^{-is}$,
thus $\cals=\cals_z$ is given by $\cals(s)=1-e^{-is}$. 
We have $f(z) = \cals(\pi) = 2$.
The curve $\gamma$ is a small circle, centered at $(2.1875,0)$ and of
radius $0.1875$, followed in the trigonometric direction. 
The snake charmer algorithm has been solved using the method 
of~\ref{odedeq2}. After one turn,
the snake slightly leans to the left. Figure A below shows the
snake after various numbers of turns of $\gamma$. One sees that
after 326 turns, the snake seems being back in its initial position.

\iffalse
\newpage
\begin{minipage}{120mm}
\hskip -5mm
\begin{minipage}[t]{40mm}
\scalebox{0.32}{
\includegraphics*[10mm,203mm][130mm,260mm]{ncer000.eps}}

\vskip -2mm\hskip 8mm Initial snake
\end{minipage}
\begin{minipage}[t]{40mm}
\scalebox{0.32}{
\includegraphics*[10mm,230mm][120mm,290mm]{ncer020.eps}}

\vskip -2mm\hskip 8mm after 20 turns
\end{minipage}
\begin{minipage}[t]{40mm}
\scalebox{0.32}{
\includegraphics*[5mm,230mm][118mm,290mm]{ncer042.eps}}

\vskip -2mm\hskip 8mm after 42 turns
\end{minipage}
\end{minipage}

\vskip 5mm
\begin{minipage}{120mm}
\hskip -5mm
\begin{minipage}[t]{40mm}
\scalebox{0.32}{
\includegraphics*[15mm,210mm][120mm,260mm]{ncer074.eps}}

\vskip -2mm\hskip 8mm after 74 turns
\end{minipage}
\begin{minipage}[t]{40mm}
\scalebox{0.32}{
\includegraphics*[15mm,210mm][120mm,260mm]{ncer162.eps}}

\vskip -2mm\hskip 8mm after 162 turns
\end{minipage}
\begin{minipage}[t]{40mm}
\scalebox{0.32}{
\includegraphics*[15mm,210mm][120mm,260mm]{ncer255.eps}}

\vskip -2mm\hskip 8mm after 255 turns
\end{minipage}

\vskip 5mm
\begin{minipage}{120mm}
\begin{minipage}[t]{40mm}
\hskip -5mm
\scalebox{0.32}{
\includegraphics*[15mm,220mm][120mm,290mm]{ncer300.eps}}

\vskip -5mm\hskip 8mm after 300 turns
\end{minipage}
\begin{minipage}[t]{40mm}
\scalebox{0.32}{
\includegraphics*[10mm,203mm][130mm,260mm]{ncer326.eps}}

\vskip -5mm\hskip 4mm after 326 turns
\end{minipage}
\end{minipage}
\end{minipage}
\fi

%\newpage

\bigskip
\noindent\hskip -30mm
\begin{minipage}{\textwidth}
%\begin{center}
\includegraphics*[10mm,153mm][230mm,260mm]{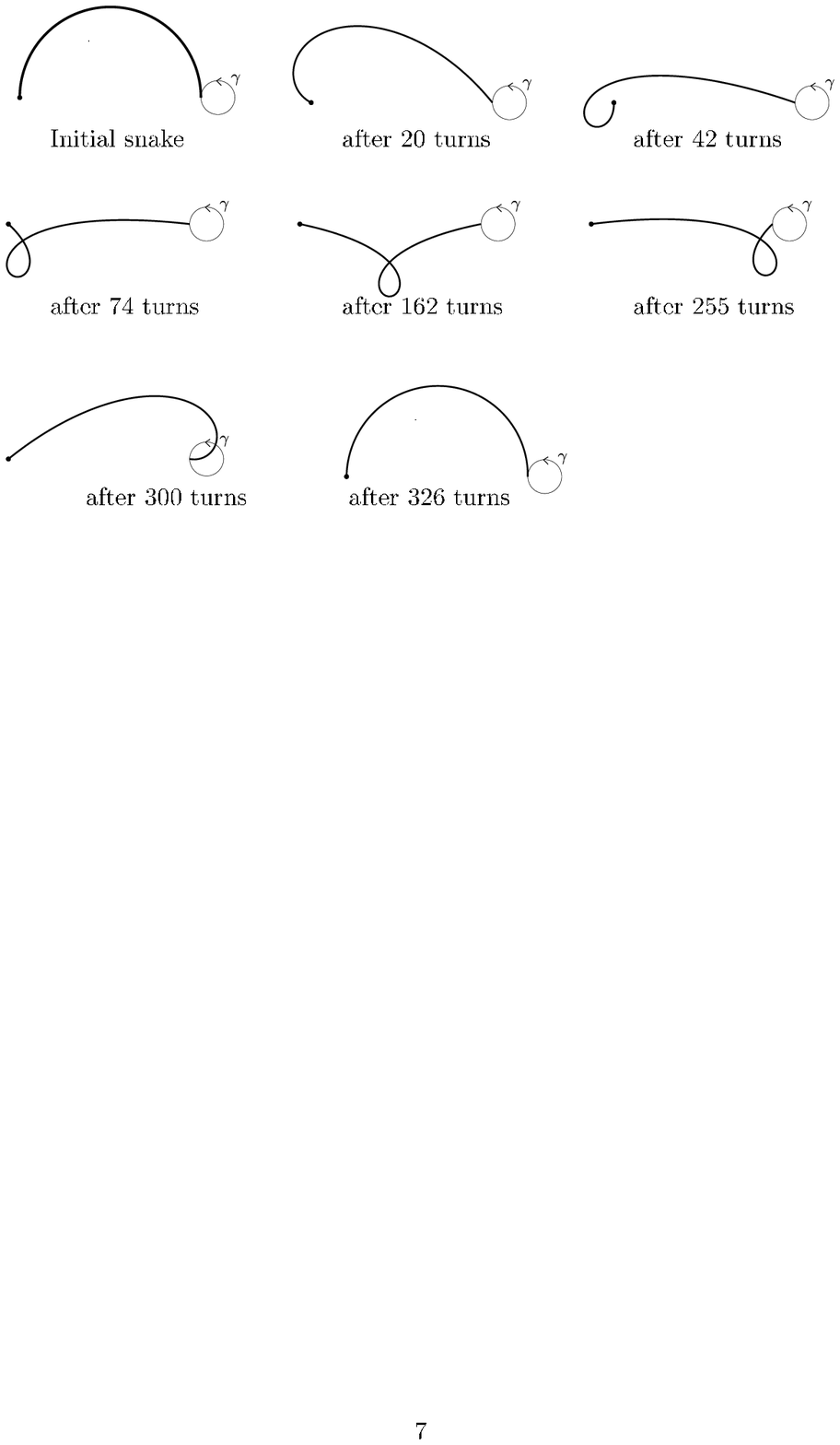}
\sk{-22}\hskip 85mm
Figure A
\end{minipage}

%\begin{center}
%\sk{-30}
%Figure A
%\end{center}
\sk{6}

The curve $t \mapsto g(t)\in\moeb{1}$ obtained from
Equation~\eqref{ode_mob} (such that $z_t=g(t)\cdot z_0$),
may be visualized using the
diffeomorphism $\bbr^2\times \bbs^1\hfl{\approx}{}\moeb{1}$
described at the end of~\ref{moeb}. The two fold covering
$SU(1,1)\to\bbr^2\times \bbs^1$ thus obtained is given
in formula by
$$
\begin{pmatrix}
a & b \\ \bar{b} & \bar{a}
\end{pmatrix}
\longmapsto
(v,e^{i\theta})
$$
with $\theta = 2\arg(a)$ and $v = 2\,\mathrm{arccosh}(|a|)e^{i\arg(ab)}$.
Figure B below illustrates the subset
$\{(v(n \cdot 2\pi),\theta(n \cdot 2\pi)) \mid n=1,\dots,326 \}
\subset\bbr^2\times\bbs^1$.
%These are the points
%of the curve $t \mapsto g(t) = (v(t),e^{\theta(t)})$ such that $g(t) \cdot z \in \ho(z)$.
The picture on the left hand side shows the 326 points of the set
$\{ v(n \cdot 2\pi) \}$
($v=(v_1,v_2)$) and the one on the right hand side is the graph of the function
$n \mapsto \theta(n \cdot 2\pi)$ (also formed by 326 points, but
hardly distinguishable).
That the points look more concentrated between 80 and 250 seems
to be related to the fact, seen in Figure A, that the snake's shape changes less
drasticly in this range. As in Figure A, we observe that $\theta(t)$ and $v(t)$ seem to have returned to their original position after 326 turns.

\bigskip
\noindent\begin{minipage}{\textwidth}
\begin{center}
\epsfig{file=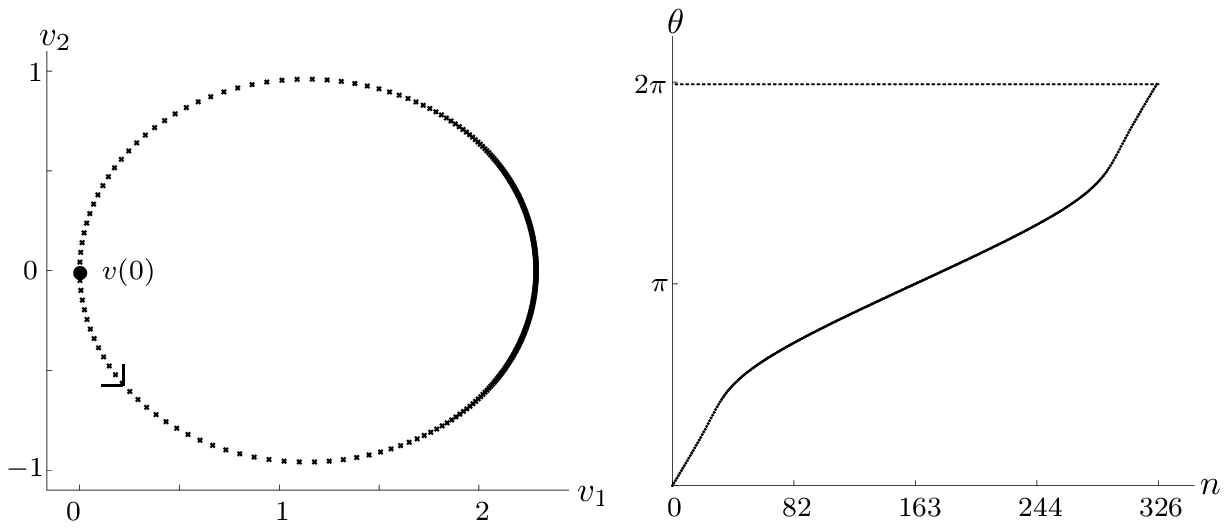,width=\textwidth}\\
Figure B
\end{center}
\end{minipage}
\bigskip

Figure C below shows the entire curve $(\theta(t),v(t))$ for $t \in [0,326\cdot 2 \pi]$.

\bigskip
\noindent\begin{minipage}{\textwidth}
\begin{center}
\epsfig{file=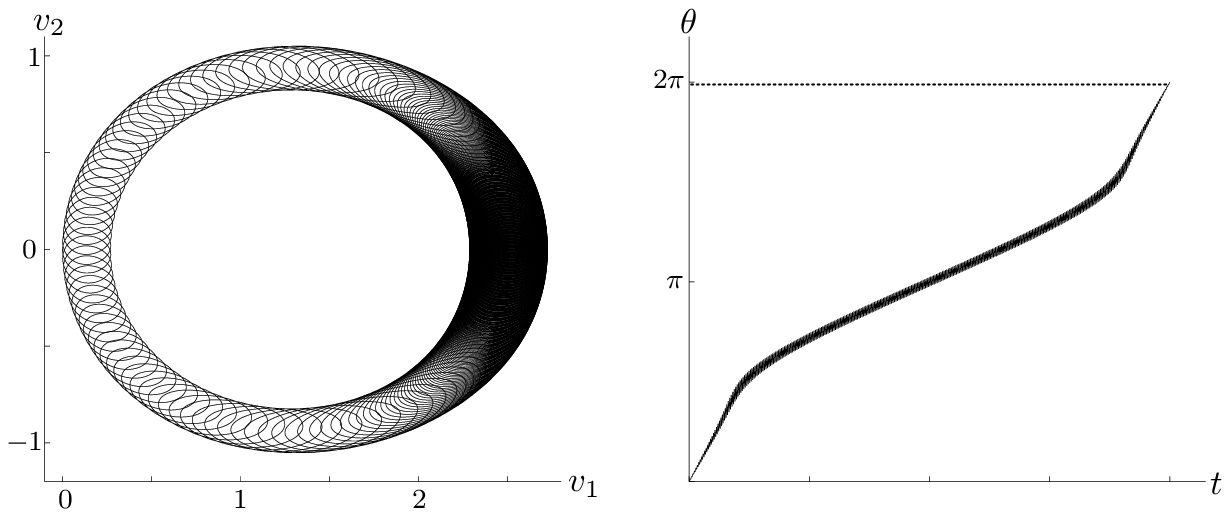,width=\textwidth}\\
Figure C
\end{center}
\end{minipage}

\end{ccote}

\begin{ccote}\label{exnoconnhorb} A non-connected holonomy orbit. \rm\
Let $z_0:[0,2\sqrt{2}]\to \bbs^1$ be the $2$-valued configuration
$$
z_0(s)=\left\{
\begin{array}{lcr}
e^{i\pi/4} & \hbox{if $s<\sqrt{2}$} \\
e^{-i\pi/4} & \hbox{if $s\geq\sqrt{2}$} \ .
\end{array}\right.
$$
Thus, $f(z)=2$. Let $\gamma:[0,1]\to\bbc$ be the
piecewise $\calc^\infty$-loop at $2$ defined by
$$
\gamma(t) = \left\{
\begin{array}{lll}
2\,e^{2i\pi t} & \hbox{if $t\leq 1/2$} \\
8(t-3/4) & \hbox{if $t\geq 1/2$} \ .
\end{array}\right.
$$

\begin{center}
%\begin{minipage}{15cm}
%\vskip -80mm
\setlength{\unitlength}{1mm}
\begin{pspicture}(5,2.5)(-6,-1)
\psline[linewidth=2.5pt](0,0)(1,1)(2,0)
\psline[linewidth=2.5pt,linecolor=lightgray](0,0)(1,-1)(2,0)

\put(-0.3,-0.4){{\large $0$}}
\put(0,0){\circle*{0.2}}
\psarc[linewidth=0.9pt,linestyle=dotted](0,0){2}{0}{180}
\psline[linewidth=0.9pt,linestyle=dotted](-2,0)(2,0)
\psline[linewidth=0.8pt](0.4,1.8)(0,2)(0.4,2.1)
\psline[linewidth=0.8pt](-1.3,0.2)(-1,0)(-1.3,-0.2)
\put(0.9,1.9){{\Large $\gamma$}}
\put(0,0.9){{\Large $\cals_{z_0}$}}
\put(0,-0.9){{\Large\gray $\cals_{z_1}$}}
\put(2.1,-0.4){{\large $2$}}
\put(-2.3,-0.4){-{\large $2$}}
\end{pspicture}
%\end{minipage}
\end{center}

By \ref{bival}, the unique possible lifting $z_t$, 
starting at $z_0$ and determined by $f(z_t)=\gamma(t)$ is horizontal.
It follows that $z_1=\bar z_0$ (the complex conjugate of $z_0$).
Hence $\ho(z_0)=\{z_0,\bar z_0\}$ is not connected.

Using an arbitrarily small bump around $0$, one can
perturb $\gamma$ so that it avoids the origin. But $z_1$ would then
be equal to $z_0$. So $z_t$ does not depend continuously on $\gamma$,
contrarily to what happens in \proref{continuite_gamma}. Here, two hypotheses
of \proref{continuite_gamma} are not satisfied: $\spdim(z_0)=0$ and,
what seems more serious, the condition $|f(z_0)| < L - 2\,\sed(z_0)$
is impossible since $\sed(z_0)=\sqrt{2}=L/2$.
\end{ccote}

\begin{ccote}\label{bivalbis}\rm  
\exref{exnoconnhorb} is an illustration of the following
computation of $\ho(z_0)$ when $z_0$ is bivalued.

\begin{Proposition}\label{ho_bival}
Let $z_0$ be bivalued configuration and let $b=f(z_0)$. Then
\renewcommand{\labelenumi}{(\roman{enumi})}
\begin{enumerate}
\item $\ho(z_0)\approx \bbs^{d-2}$ if $z_0$ is not lined;
\item $\ho(z_0)$ is one point if $z_0$ is lined and $f(z_0)\neq 0$;
\item $\ho(z_0)\approx \bbs^{d-1}$ if $z_0$ is lined and $f(z_0)=0$.
\end{enumerate}
\end{Proposition}

\preu
Let $p_0$, $q_0$, $L_p$, $L_q$, $\hat f:W\to\bbr^d$ and $W^0$ as in~\ref{bival}.
One has $\ho(z_0)\subset\hat f^\mun(b)$.
If $w\in\bbr^d$, denote by $B_w$ be the stabilizer of $w$ in $SO(d)$.
We divide the proof into four cases.

\sk{1}\noindent
{\em Case (i) with $d=2$.} \
Let $\gamma$ be a loop at $b$, with $|\gamma(t)|<L$, such that
$\gamma$ hits the sphere of radius
$|L_p-L_q|$ in a single point, tangentially and respecting
the curvature formula~\eqref{Hairereq}. By~\ref{bival},
$\gamma$ has a $\Delta$-horizontal lifting $z_t$
in $W$ and $z_t$ intersects the submanifold $W^0$
of lined configurations in one point and transversally.
This implies that $z_1\neq z_0$, which proves
Case (i) when $d=2$. \exref{exnoconnhorb} illustrates this argument
with $L_p=L_q$.

\sk{1}\noindent
{\em Case (i) with $d>2$.} \
The vectors $p$, $q$ and $b$ are co-planar and $p\neq b\neq 0$, so
$B_b\cdot z_0\approx B_b/(B_b\cap B_p)\approx SO(d-1)/SO(d-2)\approx \bbs^{d-2}$.
One has $B_b\cdot z_0\subset \ho(z_0)\subset\hat f^\mun(b)$.
Case (i) with $d>2$ then follows from the fact that $\hat f^\mun(b)\approx \bbs^{d-2}$.
The latter is proven in~\cite[Prop. 4.1]{Ha1} but, for the convenience of the reader,
we repeat the argument here in our language.
Let $W_b=\hat f^\mun(\bbr_{>0}\,b)$,
which is a submanifold of $W$ of dimension $d-1$ (see~\cite[(1.3)]{Ha1}).
Let $h:W_b\to\bbr$ defined by $h(z)=|\hat f(z)|$, which is a Morse function
whose critical point are the lined configurations of $W_b$
(see~\cite[Theorem 3.2]{Ha1}). If $z\in W_b$ is a critical point of $h$
with $h(z)\geq b$, then $z=(b/|b|,b/|b|)$, a non-degenerate maximum,
and $h(z)=L>|b|$. Therefore, $\hat f^\mun(b)=h^\mun(|b|)$ is diffeomorphic
to $\bbs^{d-2}$ by the Morse Lemma.

\sk{1}\noindent
{\em Case (ii).} \  This case is trivial since $\hat f^\mun(b)$
consists of the single point $(p_0,-p_0)$.

\sk{1}\noindent
{\em Case (iii).} \ The map $\bbs^{d-1}\to W$ given by $p\mapsto (p,-p)$
is an embedding with image $\hat f^\mun(0)$, as well as the
inclusion $B_0\cdot z_0 \subset \ho(z_0)\subset \hat f^\mun(0)$. \cqfd

\end{ccote}

%%%%%%%%%%%%%%%%%%%%%%%%%%%%%%%%%%%%%%%%%%%%%%%%%%%%%%%%%%%%%%%%%%%%%%%%%%%%%%
%\newpage

\sk{2}\noindent\small\sl
University of Geneva, Switzerland.\\
hausmann@math.unige.ch\ , eugenio.rodriguez@math.unige.ch\, .

\end{document}